\documentclass[letterpaper,11pt,oneside,reqno]{amsart}
\usepackage{amsmath,amsthm,amsfonts,amssymb,amsbsy, amscd, bbm,color}
\usepackage{graphicx,psfrag,subfigure,url}
\usepackage{cite}
\usepackage{mathrsfs}
\usepackage[colorlinks=true]{hyperref}
\usepackage[DIV13]{typearea}
\usepackage{comment}
\usepackage{caption}
\usepackage{setspace}
\usepackage{url}
\usepackage{float}
\usepackage{color}
\numberwithin{equation}{section}

\newcommand{\e}[0]{\varepsilon}

\newcommand{\R}{\ensuremath{\mathbb{R}}}

\newcommand{\Z}{\ensuremath{\mathbb{Z}}}

\newtheorem{theorem}{Theorem}[section]

\newtheorem{remark}[theorem]{Remark}

\title[Comments]{Comments on David Aldous and Persi Diaconis' ``Longest increasing subsequences: from patience sorting to the Baik-Deift-Johansson theorem''}

\author[I. Corwin]{Ivan Corwin}
\address{I. Corwin, Columbia University,
Department of Mathematics,
2990 Broadway,
New York, NY 10027, USA}
\email{ivan.corwin@gmail.com}

\begin{document}

\begin{abstract}
This is a commentary on the article: David Aldous and Persi Diaconis, \emph{Longest increasing subsequences: from patience sorting to the Baik-Deift-Johansson theorem}, {\bf Bull. Amer. Math. Soc.} 36 (1999), no. 4, 413–432. This commentary was published in {\bf Bull. Amer. Math. Soc.} 55 (2018), 363-374.
\end{abstract}

\sloppy \maketitle


\medskip

What is the distribution of the longest increasing subsequence of a permutation chosen uniformly at random from all $N!$ permutations on the numbers $\{1,\ldots, N\}$ (e.g. for a permutation $\sigma = (1,3,6,2,5,4)$ a longest increasing subsequence could be $(1,3,6)$ or $(1,2,4)$ -- both of length $3$)?
It is this question which Aldous and Diaconis took up in their 1999 AMS Bulletins article \cite{AldousDiaconis}. At the time, there had been recent breakthrough work of Baik, Deift and Johansson \cite{BDJ} which related the large $N$ asymptotic behavior of the distribution of this length to the large $N$ asymptotic behavior of eigenvalues of Gaussian Hermitian random matrices---the so called GUE Tracy-Widom distribution \cite{TracyWidomGUE} (see \cite{Forrester, AGZ} for more on random matrix theory).

Over the twenty years which have elapsed since that work, the answer to this seemingly innocent question has become intertwined with more than a few beautiful stories touching disparate areas within mathematics (e.g. probability, algebraic combinatorics, representation theory, integrable systems, algebraic geometry, partial differential equations) as well as other domains of theoretical and experimental research (e.g. physics, biology, statistics, operations research).

While it will take many books to fully expound upon the developments of these past two decades (some of the developments are contained in \cite{BaikDeiftSuidan,BorodinCorwin,BorodinGorin, BorodinPetrov,BorodinOlshanskiBook,C12,Romik}), in this short commentary I will briefly mention aspects of three such stories:
\begin{enumerate}
\item Kardar-Parisi-Zhang universality,
\item Integrable probability,
\item Stochastic PDEs.
\end{enumerate}
There are many more direction which I will not have space to discuss, including connections to classical integrable systems (e.g. Painlev\'e transcendents and Riemann Hilbert problems), quantum integrable systems, random matrix theory,  random maps, enumerative geometry, asymptotic representation theory (some discussion and further references related to these subjects can be found in \cite{Baxter,Oxfordhand,Fokasetal,Kerov,Korepietal}). In the space here provided, it will be simply impossible to cite all important and relevant works, so I must beg pardon for this as well.
\medskip

\noindent {\bf Kardar-Parisi-Zhang universality.} The longest increasing subsequence problem may seem at first to be an isolated combinatorial oddity. In truth, the behavior predicted in its large $N$ asymptotics underlies the expansive {\it Kardar-Parisi-Zhang (KPZ) universality class} \cite{KPZ} which predicts the large-scale and long-time behavior of many disparate systems (see \cite{C12,C16b,HT15,Q11,QS15} and references therein). This is facilitated by mapping the problem into a stochastic interface growth model called the {\it polynuclear growth model} \cite{KS92} as follows.

\begin{figure}[h]
\includegraphics[width=.3\textwidth]{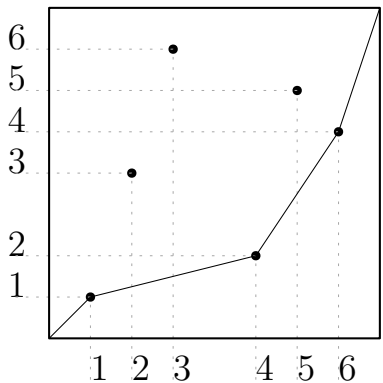}
\caption{$\mathcal{N}=6$ Poisson points in a unit square. The associated permutation is $\sigma = (1,3,6,2,5,4)$ and a longest increasing subsequence $(1,2,4)$ is drawn connecting $(0,0)$ to $(1,1)$ with lines of slope between 0 and $\pi/2$.}\label{Square}
\end{figure}

We start by {\it Poissonizing} the length of the permutation (or, in statistical physics parlance, we go to the {\it grand canonical ensemble}). Let $\mathcal{N}$ be random and {\it Poisson} distributed so that $\textrm{Prob}(\mathcal{N} = n) = e^{-N} N^n/n!$. The mean and variance of the Poisson random variable $\mathcal{N}$ are both $N$, showing that $\mathcal{N}$ concentrates  closely around $N$. Now, let $\sigma$ represent a uniform random permutation on $\mathcal{N}$ elements (first choose $\mathcal{N}$ and then choose $\sigma$ uniformly over permutations of $\{1,\ldots, \mathcal{N}\}$). There is an alternative graphical way (see Figure \ref{Square}) to sample $\sigma$. Consider a {\it Poisson point process} of intensity $N$ on a unit square $[0,1]\times [0,1]$ (this means that for each very small area $dxdy$, the probability of finding one point is approximately $Ndxdy$). The points can be ordered in terms of their $x$ coordinate and $y$ coordinate. This immediately yields a permutation (from $x$ order to $y$ order) and by construction it is equal in distribution to $\sigma$. Note that ties in this ordering happen with zero probability due to the continuous nature of the distributions.

\begin{figure}[h]
\includegraphics[width=.9\textwidth]{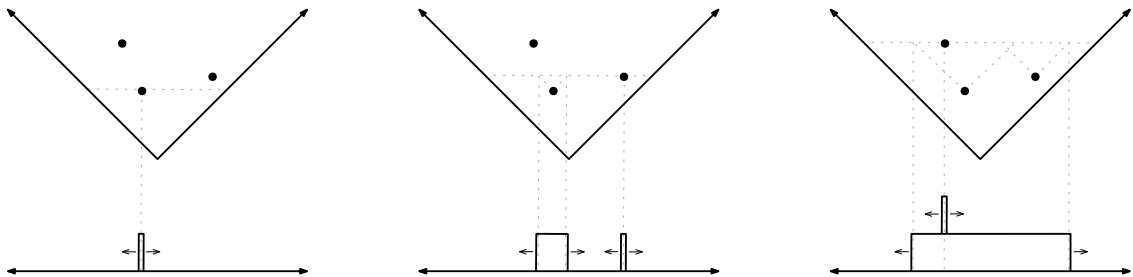}
\caption{The mapping between the graphical representation for the longest increasing subsequence and the polynuclear growth model. The three figures represent different times (as the horizontal dotted line rises vertically, time increases).}\label{LIS}
\end{figure}

The length of the longest increasing subsequence of $\sigma$ can be read off as the length of the longest path (measured by the number of Poisson points encountered along it) which starts at $(0,0)$, ends at $(1,1)$ and can only go at angles between $0$ and $\pi/2$. The beauty of this graphical representation is that it introduces a {\it space-time} dimension to the problem. Rotating the set of points by $\pi/4$ (see the top portion of Figure \ref{LIS}), one can think of the vertical direction as a {\it time} axis and the horizontal direction as a {\it space} axis.

In this rotated frame, it is natural to consider a rescaled and extended Poisson points which has intensity 1 (average number of points per unit area) and is defined on the whole cone $\{(t,x):|x|<t\}$. Let $h(t,x)$ represent the length of the longest path (where paths now need to go at angles between $\pi/4$ and $3\pi/4$). For a given $t$, $x\mapsto h(t,x)$ gives a {\it height function}  (see the bottom portion of Figure \ref{LIS}) which is zero for $|x|\geq t$, nonnegative everywhere, and increases or decreases by one at a finite number set of $x$'s.

\begin{figure}[h]
\includegraphics[width=.7\textwidth]{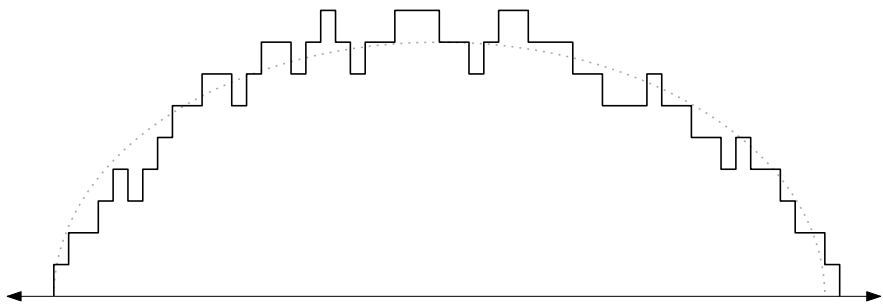}
\caption{A possible instance of the polynuclear growth model after many nucleations. The dotted line represents the circular limit shape.}\label{limit}
\end{figure}

The height function $h(t,x)$ evolves as a Markov process in time, and is called the {\it polynuclear growth model}. Specifically, a Poisson points at $(t,x)$ produces a {\it nucleation} at time $t$, position $x$ so that $h(t+,x) -h(t-,x) = 1$ ($t\pm$ represents a moment after and before $t$). The left edge and right edge of the nucleation step then move horizontally to the left and right at speed one. When a left and right moving step collide, they merge and stop moving. These dynamics are illustrated in Figure \ref{LIS}. The restriction on the Poisson points to live in $\{(t,x):|x|<t\}$ creates a growth {\it droplet} which starts at the origin and moves outwards over time.  After a long period of time and after scaling the picture back to unit size, the droplet will have grown to roughly have a circular shape (as depicted in Figure \ref{limit}).

The fluctuations of the height function grow over time and display a non-trivial spatio-temporal correlation structure. In fact, setting $h_{\epsilon}(t,x) = \epsilon h(\epsilon^{-3} t, \epsilon^{-2} x) - \bar{h}_{\epsilon}(t,x)$ for a specific centering function $\bar{h}_{\epsilon}(t,x)$ (depending on the circular limit shape), it is believed that $\lim_{\epsilon\to 0} h_{\epsilon}(\cdot,\cdot)$ exists (with the limit taking place relative to a suitable topology on functions) and is a space-time Markov process known as the KPZ fixed point \cite{CQR15,MQR17}. Moreover, it is predicted that not just this interface growth model, but a very wide class of models will all converge to this same universal fixed point (up to some model specific constants). This defines the strongest definition for membership in the {\it KPZ universality class}.

Despite great interest, the KPZ fixed point is still something of nebulous object. There are many formulas known (see \cite{MQR17} and the {\it integrable probability} discussion below) which encode distributions and transition probabilities associated to it. What is lacking for the fixed point is a simple physical characterization akin to the characterization of Brownian motion as the unique (almost surely) continuous process with independent increments and linear quadratic variation. To our knowledge, the KPZ fixed point is not the solution to a stochastic PDE, nor is it characterized via a martingale problem. Even giving an abstract set of properties which uniquely identifies the KPZ fixed point is beyond current state-of-the-art.

Instead of showing convergence to the full KPZ fixed point, a weaker notion of universality asks that a model displays the associated $3:2:1$ relationship between scaling exponents for $\textrm{time}:\textrm{space}:\textrm{fluctuation}$, or that certain marginal distributions of the height function converge to the same limiting distributions. The initial work of Baik, Deift and Johansson \cite{BDJ} identified, for the first time, analytic formulas for what one of these special universal distributions should be. Namely, they showed that the height $h(t,0)$ above the origin of the polynuclear growth model with droplet initial condition has (after correct centering) $t^{1/3}$ scale fluctuations which are asymptotically $\textrm{GUE}$ Tracy-Widom distributed. This distribution first arose in the early 90s in the context of random matrix theory \cite{TracyWidomGUE}.

There has been some progress in demonstrating the KPZ scaling exponents for more general growth models (e.g. \cite{BS10,QV07}). On the other hand, progress towards proving the distributional (e.g. $\textrm{GUE}$ Tracy-Widom) limits  has been limited to special models (like the polynuclear growth mode) which enjoy enhanced algebraic structure coming from connections to representation theory and integrable systems. These models are sometimes called {\it exactly solvable} or {\it integrable probabilistic systems} and we will discuss them further below. The solvability of these models is very precarious---most perturbations, while they should not change the universal limits, will totally destroy the algebraic structure and render existing asymptotics methods inapplicable. Despite its limitations to these special models, integrable probability provides a mathematical laboratory for discovering new asymptotic phenomena, crafting precise predictions and testing numerical methods and other non-rigorous methods from physics.

Through a combination of rigorous mathematical results from integrable probability, theoretical physics methods, numerical simulations and experimental demonstrations, the notion and scope of KPZ universality has been refined and expanded quite drastically over the past decades (see the reviews \cite{C12,C16b,HT15,Q11,QS15}). In fact, through various mappings (some quite evident like the relation between the longest increasing subsequence and the polynuclear growth model, and others hidden in algebraic structure like the relation of those models to random matrix theory), the KPZ universality class has steadily grown to encompass much more than just random growth models. The universality class connects to and predicts certain asymptotic behaviors for models of traffic, queues in series, mass transport in random media, turbulence, stochastic optimization, quantum entanglement; the statistics and scaling exponents are closely related to those which arise in random matrix theory and random tilings; and experimental / numerical work has revealed these behaviors in many other physical systems like superconductor vortices, disordered liquid crystals, bacterial colony growth and competition interfaces, cancer growth, chemical reaction fronts, slow combustion, coffee stains, and conductance fluctuations in Anderson localization.

\medskip

\noindent {\bf Integrable probability.}
The longest increasing subsequence (in its various disguises through mappings described above) is an {\it integrable probabilistic system} in that (1) there are many concise and exact formulas known for expectations of observables of interest (e.g. distribution of the length), and  (2) taking a limit of the system along with the expectations and observables provides a precise description of a much wider universality class (i.e., the KPZ class).  Here `integrable' should not be confused with the notion of an integrable function (whose integral is finite). Rather, it is being used in the spirit of classical or quantum `integrable systems' in which there are many conserved quantities which allow one to `integrate' or solve the system. Another synonym, `exactly solvable', is often used in the study of 2D lattice models \cite{Baxter}---a subject closely related to this one.

\begin{figure}[h]
\includegraphics[width=.35\textwidth]{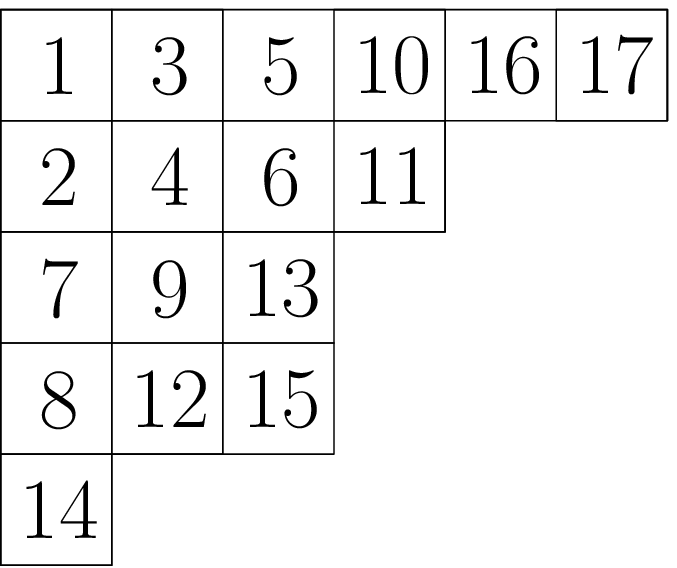}
\caption{A standard Young diagram of size $17$. Notice that every row and every column is increasing and all numbers $\{1,\ldots, 17\}$ are used only once.}\label{Young}
\end{figure}

One manifestation of the integrability of the longest increasing subsequence comes from the following very remarkable formula. Let $L_N$ denote the length of the longest increasing subsequence of permutation drawn uniformly at random from all $N!$ perumtations on $\{1,\ldots, N\}$. Then
$$
\mathrm{Prob}\big(L_N = \ell\big) =  \sum_{\lambda \vdash N: \lambda_1=\ell}  \frac{\mathrm{dim}_{\lambda}^2}{N!}
$$
where $\lambda= (\lambda_1\geq \lambda_2\geq \cdots \geq 0)$ with $\lambda_i\in \Z$ and $\sum \lambda_i =N$ is called a {\it Young diagram} or {\it partition} of size $N$, and   $\mathrm{dim}_\lambda$ is the dimension of the irreducible representation of the symmetric group indexed by the partition $\lambda$ . There is a simple combinatorial formula for $\mathrm{dim}_{\lambda}$ as the number of standard Young {\it tableaux} of shape $\lambda$ (i.e., a filling of the diagram in Figure \ref{Young} with the numbers $\{1,\ldots, N\}$ so that within in row and column, the numbers are strictly increasing). The measure $\mathrm{dim}_{\lambda}^2 / N!$ on partition $\lambda$ of size $N$ is called the Plancherel measure---see \cite{BorodinGorin,BorodinOlshanskiBook,Kerov} for more background and discussion.

The Plancherel measure enjoys many amazing properties which enable precise asymptotic calculations. Moreover it fits into a far-reaching hierarchy of measures which have many varied applications across fields  of mathematics. For instance, the Plancherel measure on partitions (in fact its Poissonization where $N$ is chosen as a Poisson random variable) is a special case of the Schur measure on partitions $\lambda$ (introduced in \cite{Okounkov} fairly soon after the work of Baik-Deift-Johansson) given by the following prescription
$$
\mathrm{Prob}\big(\lambda\big) = s_{\lambda}(a_1,\ldots, a_m)s_{\lambda}(b_1,\ldots, b_n) \prod_{i,j} (1-a_ib_j).
$$
Here $s_{\lambda}$ is the Schur symmetric polynomial (for a general reference for the symmetric functions discussed below, see \cite{Mac}). The measure depends on the $a$ and $b$ parameters at which the Schur polynomials are evaluated, and the fact that it sums to one over all partition $\lambda$ is the content of the classical Cauchy-Littlewood identity for Schur polynomials. Due to the Schur polynomial branching rule, the Schur measure can be related to the study of random tiling models, and through certain limit transition it relates to the Gaussian Unitary Ensemble (GUE) measure on eigenvalues in random matrix theory. The Schur measure is a {\it determinantal point process} (as can be showed from determinant and non-crossing path formulas for Schur polynomials) which provides a well-developed calculus for computing asymptotics. In the decade after the Baik-Deift-Johansson result, Schur measure and the methods developed to study it became a major tool in uncovering the behavior of the KPZ universality class---see for instance \cite{BaikDeiftSuidan,BF14, BorodinGorin,BorodinPetrov,BorodinOlshanskiBook,Romik}.

\begin{figure}[h]
\includegraphics[width=.9\textwidth]{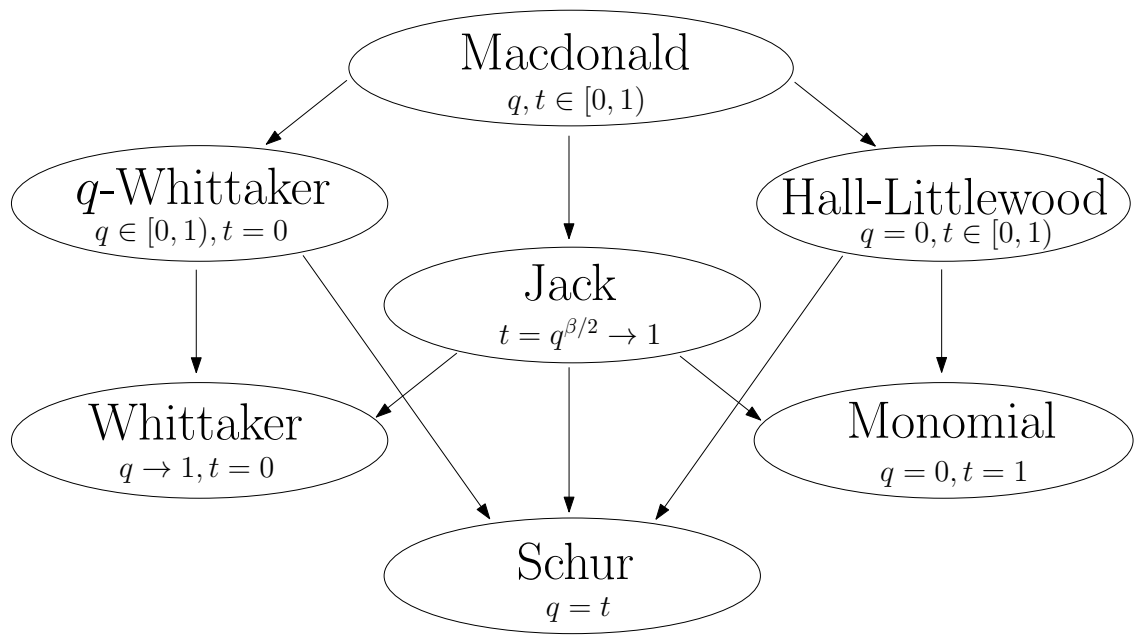}
\caption{Hierarchy of limits of Macdonald symmetric polynomials}\label{Macdonald}
\end{figure}

The hierarchy of Macdonald symmetric polynomials (see Figure \ref{Macdonald} as well as \cite{Mac}) provides a natural generalization of Schur symmetric polynomials, while maintaining or deforming many of their key properties. Macdonald polynomials involve two extra parameters $q$ and $t$ and degenerate for Schur polynomials when $q=t$. Other well studied symmetric functions which fit into the hierarchy (through various specialization of parameters or limit transitions) are Hall-Littlewood polynomials (related to representation theory over finite or $p$-adic fields), Jack polynomials (related to spherical functions for Riemannian symmetric spaces), and Whittaker functions (related to the quantum Toda lattice). The natural deformation of Schur measure is to replace the Schur polynomials by Macdonald polynomials (or other functions in the hierarchy). The study of such deformations of Schur measure, as well as the development of methods to connect these measures to interesting probabilistic systems and methods to analyze asymptotics (since the determinantal structure of Schur measure is lost in these deformations) has been a focus of significant research in the past decade---for instance \cite{BorodinCorwin,BorodinGorin, BorodinPetrov,C16b,OCon15} and references therein.

There are now many interesting probabilistic systems which fit into this framework. Schur measure relates to random plane partitions, tilings, longest increasing subsequences, polynuclear growth, TASEP, last passage percolation and GUE random matrix theory. Measures defined via monomial symmetric polynomials relate to Kingman partition structures, and cycles structures in random permutations. Jack measure relates to general $\beta$ random matrix theory ensembles. Hall-Littlewood measure relates to random matrices over finite fields (and associated Cohen-Lenstra heuristics in number theory) as well as stochastic vertex models coming from quantum integrable systems (i.e. the six vertex model and its generalizations). $q$-Whittaker measure relates to traffic models like $q$-TASEP, and Whittaker measure relates to directed polymer models (which generalize the graphical construction of the longest increasing subsequence).

Among the systems which fit into the Macdonald hierarchy, one which has received particular attention is the Kardar-Parisi-Zhang (KPZ) stochastic partial differential equation (SPDE) (or KPZ equation for short), which arises as a limit of both the Whittaker and Hall-Littlewood measures. The machinery developed to study these measures has provided a route (e.g., see \cite{BCF}) to study the exact distribution of the solution to this SPDE and prove various long-time limit theorems (analogous to that of Baik-Deift-Johansson for the longest increasing subsequence). There are other approaches (such as using tools from Bethe ansatz and quantum integrable systems) which have provided related but different routes to study the KPZ equation exact statistics (and many other probabilistic systems)---e.g., see \cite{BP15,CP14} and references therein.

\medskip

\noindent{\bf Stochastic PDEs.}
The KPZ equation describes the evolution of a height function $h(t,x)$ where $t\geq 0$ denotes time, $x\in \R$ space, and $h(t,x)\in \R$ the height above $x$ at time $t$. The SPDE is
\begin{equation}\label{KPZ}
\partial_t h(t,x) = \tfrac{1}{2}\partial_{x}^2 h(t,x) +  \tfrac{1}{2}\big(\partial_x h(t,x)\big)^2 + \xi(t,x)
\end{equation}
where $\xi(t,x)$ denotes a space-time white noise (i.e., a generalized Gaussian process which has covariance which is a delta function is space and time). If $\xi$ were replaced by a deterministic and smooth function, then the Hopf-Cole transform would relate this equation of Hamilton-Jacobi type to a heat equation so that
\begin{equation}\label{SHE}
h(t,x) = \log z(t,x)\qquad \textrm{where}\qquad \partial_z(t,x) = \tfrac{1}{2} \partial_{x}^2 z(t,x) + \xi(t,x) z(t,x).
\end{equation}
For white noise $\xi$, the above SPDE is called the stochastic heat equation (SHE). See \cite{C12,Q11} for further background and discussion related to the below content.

There are two different things one could mean by the phrase `solving the KPZ equation'. The first notion is in the spirit of PDE theory in which `solving' corresponds to constructing (generally through a fixed-point scheme) a mapping from initial data and noise onto solutions. The key questions here usually involves understanding the local properties (e.g. regularity) of such mappings, as well as determining whether global solutions exist. The second notion is in the spirit of integrable systems in which one seeks to compute explicit formulas for the distribution of the constructed solution. From the directions of integrable probability and KPZ universality, there was a great deal of attention in the last decade ago devoted to this second notion for the KPZ equation---namely, the question of determining the distribution of the solution, for example, at a given space-time point and with given initial data.

Likewise, in the past decade there has been a flurry of activity related to the PDE side (first notion above) of studying the KPZ equation, as well as many other SPDEs. Making sense of the KPZ equation \eqref{KPZ} is not particularly easy and that challenge has prompted a number of important developments. The roughness of white noise suggests that $h$ should be as rough as Brownian motion in space (e.g. if one removed the non-linearity, then the linear equation is easy to solve and one sees this holds true). This roughness poses a problem since the non-linearity asks to differentiate and then square the rough process. So, one has to work to overcome this challenge.

The simplest scheme to define the solution is to define the {\it Hopf-Cole solution} to the KPZ equation as in \eqref{SHE}. In other words, solve the SHE with initial data $e^{h(0,x)}$ and then take the logarithm of the solution $z(t,x)$ to get $h(t,x)$. This turns out to be the physically relevant notion of solution to the KPZ equation---it captures all of the expected behavior and arises quite universally from various approximation schemes to the equation. While the problem of showing that the KPZ equation arises from all sorts of smoothed or discretized systems has sparked a number of developments recently, the first such results go back to the early 90s. It was first demonstrated in \cite{BertiniCancrini} that the KPZ equation with a smoothed noise (given by convolving $\xi$ against a bump function in space) minus an appropriate function (depending on time and the smoothing) converges as the smoothing disappears to the Hopf-Cole solution. The key here is that once the noise is smoothed, one can make direct sense of the equation.

\begin{figure}[h]
\includegraphics[width=.4\textwidth]{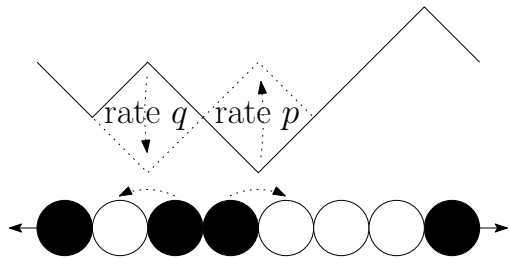}
\caption{The ASEP height function $h$ (on top) and associated particle occupation function $\eta$ (on bottom). Height function local maxima turn to local minima at rate $q$ and the opposite at rate $p$; particles hop left at rate $q$ and right at rate $p$.}\label{asep}
\end{figure}

The first discrete system which was proved to converge to the KPZ equation was the {\it asymmetric simple exclusion process} (ASEP) which is a simple growth model in which the height function $h^{\textrm{ASEP}}(t,x)$ takes integer values for $x\in \Z$, and satisfies $\big|h^{\textrm{ASEP}}(t,x) - h^{\textrm{ASEP}}(t,x+1) \big|=1$ for all $t\geq 0$ and $x\in \Z$. The function evolves by changing local minimums in the height function into local maximums at rate (i.e., after exponential random weighting times of rate) $p$ and the opposite at rate $q=1-p$. See Figure \ref{asep} for an illustration.

In order for ASEP to converge to the KPZ equation it is necessary to introduce {\it weak asymmetry} scaling. For the sake of stating a result, consider ASEP started so that $h^{\textrm{ASEP}}(0,0)=0$ and $\big\{h^{\textrm{ASEP}}(0,x) - h^{\textrm{ASEP}}(0,x+1)\big\}_{x\in \Z}$ is distributed as a collection of independent $\pm 1$ random variables with equal probability to both outcome (i.e. $h^{\textrm{ASEP}}(0,x)$ is a random walk in $x$). Now, consider a sequence of ASEP models where the asymmetry $p-q=\e$ and the height function are rescaled with $\e$ as $h^{\textrm{ASEP}}_{\e}(t,x) = \e h^{\textrm{ASEP}}(\e^{-4} t,\e^{-2} x) - \e^{-2} t/2$. Then, \cite{BertiniGiacomin} proved that $h^{\textrm{ASEP}}_{\e}$ converges as $\e\searrow 0$ to the Hopf-Cole solution to the KPZ equation (i.e., the measure on space-time functions induced by $h^{\textrm{ASEP}}_{\e}$ converges in distribution to that induced by $h$). It should be noted that if one did not perform the special weak asymmetry, then ASEP is believed to converge to the KPZ fixed point, not the KPZ equation. This was been demonstrated in \cite{TWASEP} at the level of one-point marginal distributions using tools from integrable probability.

The convergence result of \cite{BertiniGiacomin} relied on the fact that the ASEP height function satisfies a discrete version of the Hopf-Cole transform, under which its microscopic dynamics transform into a microscopic version of the SHE. Since the SHE is linear, it is then much easier to prove convergence at this level. The existence of a microscopic Hopf-Cole transform is a delicate matter---changing the dyamics of ASEP to allow for more general growth rules or height changes will generally break this mechanism. Integrable probability provides a systematic way to find discrete systems which enjoy similar exact transforms. Indeed, the microscopic Hopf-Cole transform is a consequence of something called {\it Markov duality} whose origins can be traced back to symmetries of quantum groups \cite{Schutz,CGRS}. Thus, the Hopf-Cole transform approach has worked for the handful of systems enjoying such dualities, and effectively proved their convergence to the KPZ equation.

Define  $\eta(t,x) = h^{\textrm{ASEP}}(t,x) - h^{\textrm{ASEP}}(t,x+1)$ and think of $\eta(t,x)=1$ as corresponding to there being a particle at $x$ and $\eta(t,x)=-1$ as there being a hole. Then ASEP can be interpreted as a simple particle hopping model. Particles attempt to jump left at rate $p$ and right at rate $q$, but only make the jump is the destination is unoccupied. From this perspective, it is clear that one can generalize the model considerably---for instance allowing for non-nearest-neighbor jumps, or jump rates which depend on local configurations. The Hopf-Cole approach is, for the most part, too rigid to accommodate these changes (though \cite{Dembo2016} managed to extend the approach to some non-nearest-neighbor systems).

There are now a number of other approaches to show the robustness of the KPZ equation under various approximation schemes. Besides the Hopf-Cole method, these include the methods of energy solutions \cite{GJ10,GP2015a}, regularity structures \cite{Hairer13,Hairer14}, paracontrolled distributions \cite{GP15b,GIP}, and renormalization group \cite{Antti}. Each approach offers some advantages and some limitations, though they all yield the same Hopf-Cole solution. For instance, energy solutions can be used to prove the general non-nearest-neighbor exclusion processes converge (under weak asymmetry) to the KPZ equation (see, e.g. \cite{GJS15})---but the convergence can only be done for processes started from their stationary measure (i.e., where the $\eta(0,x)$ are $\pm 1$ independently in $x$ with the same probabilities). The other three notions have been useful in showing that stochastic Hamilton-Jacobi equations (i.e., replacing the $\big(\partial_x h\big)^2$ terms with a general non-linear function) or equations with non-Gaussian noise all converge to the KPZ equation when rescaled appropriates (see, e.g. \cite{HaiQua,HaiShen}). Those approaches and results are, however, presently limited to periodic domains and cannot be used to prove convergence results on the full line. Presently, the scope of each of these approaches has non-trivial symmetric difference---each approach offers something new and useful to understanding the approximation theory for the KPZ equation. Moreover, some of these approaches (namely regularity structures, paracontrolled distributions and renormalization group) can be applied systematically for more general classes of non-linear SPDEs.

\medskip

\noindent {\bf Conclusion.}
Twenty years ago, Baik, Deift and Johansson discovered a remarkable limit theorem connecting asymptotics of the longest increasing subsequence to those of random matrix theory. Aldous and Diaconis related a few stories (e.g. patience sorting) related to this theorem.  In this commentary, I have tried to relate few other stories that developed over the elapsed time in a manner closely intertwined with this original work (but not covered in Aldous and Diaconis' article). The number of topics and fields I have touched upon is far exceeded by those for which I have said nothing. Surely, the interconnectedness of all of these areas is a beautiful illustration of the unity of mathematics and science.

\medskip
\noindent {\it Acknowledgements}
The author wishes to thank Herbert Spohn for helpful remarks on a draft of this commentary. The author recognizes partial NSF support DMS-1664650, as well as support from a Packard Foundation Fellowship for Science and Engineering.


\begin{thebibliography}{alpha}

\bibitem{AldousDiaconis}
D. Aldous and P. Diaconis.
\newblock Longest increasing subsequences: from patience sorting to the Baik-Deift-Johansson theorem.
\newblock {\it Bull. Amer. Math. Soc.}, {\bf 36}: 413--432, 1999.


\bibitem{AGZ}
G. W. Anderson, A. Guionnet and O. Zeitouni.
\newblock {\it An Introduction to Random Matrices}.
\newblock Cambridge University Press, 2009.



\bibitem{BDJ}
J. Baik, P. Deift and K. Johansson.
\newblock On the distribution of the length of the longest increasing subsequence of random permutations.
\newblock {\it J. Amer. Math. Soc.}, {\bf 12}: 1119--1178, 1999.

\bibitem{BaikDeiftSuidan}
J. Baik, P. Deift and T. Suidan.
\newblock {\it Combinatorics and Random Matrix Theory}.
\newblock AMS Graduate Studies in Mathematics {\bf 172}, 2016.


\bibitem{BS10}
M. Bal\'azs and Timo Sepp\"{a}l\"{a}inen.
\newblock Order of current variance and diffusivity in the asymmetric simple exclusion process.
\newblock {\it Ann. Math.}, {\bf 171}:1237--1265, 2010.


\bibitem{Baxter}
R. Baxter.
\newblock {\it Exactly solved models in statistical mechanics}.
\newblock Dover, 2008.


\bibitem{BertiniCancrini}
L.~Bertini and  N.~Cancrini.
\newblock  The stochastic heat equation: Feynman-Kac formula and intermittence.
\newblock {\it J. Stat. Phys.}, {\bf 78}:1377--1401, 1995.

\bibitem{BertiniGiacomin}
L.~Bertini and G.~Giacomin.
\newblock Stochastic {B}urgers and {KPZ} equations from particle systems.
\newblock {\it Commun. Math. Phys.}, {\bf 183}:571--607, 1997.


\bibitem{BorodinCorwin}
A.~Borodin and I.~Corwin.
\newblock Macdonald processes.
\newblock {\it Probab. Theor. Rel. Fields}, {\bf 158}:225--400, 2014.

\bibitem{BCF}
A.~Borodin, I.~Corwin, P.~L.~Ferrari.
\newblock Free energy fluctuations for directed polymers in random media in $1+1$ dimension.
\newblock {\it Commun. Pure Appl. Math.}, {\bf 67}:1129--1214, 2014.


\bibitem{BF14}
A. Borodin and P. L. Ferrari.
\newblock Anisotropic growth of random surfaces in $2+1$ dimensions.
\newblock {\it Commun. Math. Phys.}, {\bf 325}:603--684, 2014.


\bibitem{BorodinGorin}
A. Borodin and V. Gorin,
\newblock Lectures on integrable probability.
\newblock In: {\it Probability and Statistical Physics in St. Petersburg, Proceedings of Symposia in Pure Mathematics}, {\bf 91}:155–214, 2016.

\bibitem{BorodinPetrov}
A. Borodin and L. Petrov,
\newblock Integrable probability: From representation theory to Macdonald processes
\newblock {\it Probab. Surveys}, {\bf 11}:1--58, 2014.


\bibitem{BP15}
A. Borodin and L. Petrov.
\newblock Integrable probability: stochastic vertex models and symmetric functions
\newblock In {\it Stochastic Processes and Random Matrices: Lecture Notes of the Les Houches Summer School} {\bf 104}, 2015.


\bibitem{BorodinOlshanskiBook}
A. Borodin and G. Olshanski.
\newblock {\it Representations of the Infinite Symmetric Group}.
\newblock Cambridge Studies in Advanced Mathematics, 2016.


\bibitem{CGRS}
G.~Carinci, C.~Giardin\`a, F.~Redig and T.~Sasamoto.
\newblock A generalized asymmetric exclusion process with {$U_q(\mathfrak{sl}_2)$} stochastic duality.
\newblock {\it Probab. Theo. Rel. Fields}, {\bf 166}:887--933, 2016.



\bibitem{C12}
I. Corwin.
\newblock The {K}ardar-{P}arisi-{Z}hang equation and universality class.
\newblock {\it Random Matrices Theory Appl.}, {\bf 1}, 2012.


\bibitem{C16b}
I.~Corwin.
\newblock Kardar-Parisi-Zhang universality.
\newblock {\it Notices of the American Mathematical Society}, March, 2016.


\bibitem{CP14}
I. Corwin and L. Petrov.
\newblock Stochastic higher spin vertex models on the line.
\newblock {\it Commun. Math. Phys.}, {\bf 343}:651--700, 2016.

\bibitem{CQR15}
I.~Corwin, J.~Quastel, D.~Remenik.
\newblock The renormalization fixed point of the Kardar-Parisi-Zhang universality class.
\newblock {\it J. Stat. Phys.}, {\bf 160}:815--834 (2015).

\bibitem{Dembo2016}
A.~Dembo and L.-C. Tsai.
\newblock Weakly asymmetric non-simple exclusion process and the
  {K}ardar--{P}arisi--{Z}hang equation.
\newblock {\it Commun. Math. Phys.}, {\bf 341}:219--261, 2016.

\bibitem{Oxfordhand}
P. Di Francesco, J. Baik, G. Akemann (editors).
\newblock {\it The Oxford Handbook of Random Matrix Theory}.
\newblock Oxford University Press, 2011.


\bibitem{Fokasetal}
S. Fokas, R. Its, A. Kapaev and V. Y. Novokshenov.
\newblock {\it Painlev\'{e} Transcendents: The Riemann-Hilbert Approach}.
\newblock AMS Mathematical Surveys and Monographs, {\bf 128}, 2006.


\bibitem{Forrester}
P. Forrester.
\newblock {\it Log-Gases and Random Matrices}.
\newblock London Mathematical Society Monographs {\bf 34}, 2010.


\bibitem{GJ10}
P.~{Goncalves} and M.~{Jara}.
\newblock {Universality of KPZ} equation.
\newblock {\it arXiv:1003.4478}, 2010.


\bibitem{GJS15}
P.~Gon{\c{c}}alves, M.~Jara and S.~Sethuraman.
\newblock A stochastic burgers equation from a class of microscopic
  interactions.
\newblock {\it Ann. Probab.}, {\bf 43}:286--338, 2015.




\bibitem{GIP}
M.~Gubinelli, P.~Imkeller and N.~Perkowski.
\newblock Paracontrolled distributions and singular {PDE}s.
\newblock {\it Forum Math. PI}, 3(1), 2015.


\bibitem{GP15b}
M.~Gubinelli and N.~Perkowski.
\newblock {KPZ} reloaded.
\newblock {\it Commun. Math. Phys.}, 349:165--269, 2017.


\bibitem{GP2015a}
M.~Gubinelli and N.~Perkowski.
\newblock {E}nergy solutions of {KPZ} are unique.
\newblock {\it J. Amer. Math. Soc.}, 31:427--471, 2017.


\bibitem{Hairer13}
M.~Hairer.
\newblock Solving the {KPZ} equation.
\newblock {\it Ann. Math.}, 178(2):559--664, 2013.

\bibitem{Hairer14}
M.~Hairer.
\newblock A theory of regularity structures.
\newblock {\it Invent. Math.}, 198(2):269--504, 2014.


\bibitem{HaiQua}
M.~Hairer and J.~Quastel.
\newblock A class of growth models rescaling to {KPZ}.
\newblock {\it arXiv:1512.07845}, 2015.

\bibitem{HaiShen}
M.~Hairer and H.~Shen.
\newblock A central limit theorem for the {KPZ} equation.
\newblock {\it Ann. Probab.}, {\bf 45}:4167--4221, 2017.



\bibitem{HT15}
T. Halpin-Healy, K. Takeuchi.
\newblock A {KPZ} cocktail -- shaken, not stirred...
\newblock {\it J. Stat. Phys.}, {\bf 160}:794--814 (2015).


\bibitem{KPZ}
K.~Kardar, G.~Parisi, Y.Z.~Zhang.
\newblock  Dynamic scaling of growing interfaces.
\newblock {\it Phys. Rev. Lett.}, {\bf 56}:889--892, 1986.

\bibitem{Kerov}
S. V. Kerov.
\newblock {\it Asymptotic Representation Theory of the Symmetric Group and its Applications in Analysis}.
\newblock AMS Translations of Mathematical Monographs, {\bf  219}, 2003.


\bibitem{Korepietal}
V. E. Korepin, N. M. Bogoliubov and A. G. Izergin.
\newblock {\it Quantum Inverse Scattering Method and Correlation Functions}.
\newblock Cambridge University Press, 2010.


\bibitem{KS92}
J. Krug and H. Spohn.
\newblock Kinetic Roughening of Growing Interfaces.
\newblock In {\it Solids far from Equilibrium: Growth, Morphology and Defects} 479--582, Cambridge University Press, 1992.

\bibitem{Antti}
A.~Kupiainen.
\newblock Renormalization group and stochastic {PDE}s.
\newblock {\it Ann. H. Poincare}, 17:497--535, 2016.




\bibitem{Mac}
I. G. Macdonald.
\newblock {\it Symmetric functions and Hall polynomials}.
\newblock Oxford Mathematical Monographs, 1995.

\bibitem{MQR17}
K. Matetski, J. Quastel and D. Remenik.
\newblock The KPZ fixed point.
\newblock arXiv:1701.00018, 2017.

\bibitem{OCon15}
N. O'Connell.
\newblock Whittaker functions and related stochastic processes.
\newblock In: {\it Random matrix theory, interacting particle systems, and integrable systems.} 385--409 Cambridge University Press, 2015.


\bibitem{Okounkov}
A. Okounkov
\newblock Infinite wedge and random partitions.
\newblock {\it Sel. Math.}, {\bf 7}:57--81, 2001.


\bibitem{Q11}
J. Quastel
\newblock Introduction to {KPZ}.
\newblock {\it Curr. Dev. Math.}, {\bf 1}, 2011.

\bibitem{QS15}
J. Quastel, H. Spohn.
\newblock The one-dimensional {KPZ} equation and its universality class.
\newblock {\it J. Stat. Phys.}, {\bf 160}:965--984, 2015.

\bibitem{QV07}
J. Quastel and B. Valk\'o.
\newblock $1/3$ Superdiffusivity of Finite-Range Asymmetric Exclusion Processes on $Z$.
\newblock {\it Commun. Math. Phys.}, {\bf 273}:379--394, 2007.

\bibitem{Romik}
D. Romik.
\newblock {\it The Surprising Mathematics of Longest Increasing Subsequences}.
\newblock Cambridge University Press, 2015.

\bibitem{Schutz}
G.~M. Sch\"utz.
\newblock Duality relations for asymmetric exclusion processes.
\newblock {\it J. Stat. Phys.}, {\bf 86}:1265--1287, 1997.


\bibitem{TracyWidomGUE}
C. A. Tracy and H. Widom.
\newblock Level-spacing distributions and the Airy kernel.
\newblock {\it Commun. Math. Phys.}, {\bf 159}: 151--174, 1994.


\bibitem{TWASEP}
C. Tracy and H. Widom.
\newblock Asymptotics in ASEP with Step Initial Condition.
\newblock {\it Commun. Math. Phys.}, {\bf 290}:129--154, 2009.







%
%
%
%
%
%
%
%
%
%
%
%
%
%
%
%
%
%
%
%
%
%
%
%
%
%
%
%
%
%
%
%
%
%
%
%
%
%
%
%

\end{thebibliography}
\end{document}